\documentclass[11pt]{article}

\usepackage{amssymb}
\usepackage{amsmath,amsthm}
\usepackage{amscd}
\usepackage{amstext}
\usepackage{amsfonts}
\usepackage[all]{xy}
\usepackage{titlesec}

\newtheorem{lemma}{Lemma}[section]
\newtheorem{theorem}[lemma]{Theorem}
\newtheorem{corollary}[lemma]{Corollary}
\newtheorem{proposition}[lemma]{Proposition}
\newtheorem{remark}[lemma]{Remark}

\theoremstyle{definition}
\newtheorem*{definition}{Definition}

\def\Ext{\mathop{\rm Ext}\nolimits}

\def\Hom{\mathop{\rm Hom}\nolimits}

\newcommand{\ci}{\cal{I}\ }
\newcommand{\cip}{\cal{I}^{\perp}\ }
\title{\Large \bf  Precovering and preenveloping Ideals
}

\author{Furuzan Ozbek\\
{Department of Mathematics, University
of Kentucky, Lexington, Kentucky, USA} }
\date{}
\begin{document}
\baselineskip=18pt \maketitle
\begin{abstract}
In  \cite{Salce} L. Salce introduced the notion of a cotorsion pair $(\cal{F},\cal{C}) $ in the category of abelian groups. But his definitions and basic results carry over to more general abelian categories and have proven useful in a variety of settings. A significant result of cotorsion theory  proven by Eklof \& Trlifaj is that if a pair $(\cal{F},\cal{C})$ of classes of $R$-modules is cogenerated by a set, then it is complete \cite{Eklof}. Recently Herzog, Fu, Asensio  and Torrecillas developed the ideal approximation theory \cite{Herzog}, \cite{idealapp}. In this article we look at a result motivated by the Eklof \& Trlifaj argument for an ideal $\ci$ when it is generated by a set of homomorphisms.

\end{abstract}

\vspace{0.2cm}
\titlelabel{\thetitle.\quad}
\section{Introduction}

The concepts of preenvelope and precover were introduced by Enochs \cite{pre} for classes of modules. Since then the definition has been applied to different classes of categories. One of the recent application is given by Herzog, Fu, Asensio and Torecillas \cite{idealapp}. Herzog first defined and looked at phantom morphisms \cite{phantom}. Some properties of phantom morphisms raise interest in the subfunctors of the bifunctor $\Hom$ which they call ideals \cite{Herzog}.  Herzog then applied the definition of preenvelope(precover), special preenvelope(special precover) and cotorsion pair to the ideal case in ~\cite{Herzog}.

We see that several nice results in category of modules such as Salce's lemma \cite{Salce} carry over to the ideal case in Herzog's paper \cite{Herzog}. One significant result in cotorsion theory proven by Eklof \& Trlifaj is that if a pair $(\cal{F},\cal{C})$ of classes of $R$-modules is cogenerated by a set, then it is complete \cite{Eklof}. We look at how this partially carries over to the ideal case.

Throughout this paper, we will focus on the ideals which are generated by a set of homomorphisms. First we examine how the elements of such ideals look like (Remark~\ref{setgenerated}). This helps us to characterize the elements of $\cip$ (Lemma~\ref{iperp}). We see that every homomorphism $g$ in $\cip$ has a small enough factorization (Lemma~\ref{factor}) which is motivated by the Proposition(5.2.2) from Enochs and Jenda's Relative Homological Algebra \cite{RHA}. With all the tools in hand we prove that if an ideal $\ci$ is generated by a set then $\cip$ is a preenveloping class (Theorem~\ref{preenveloping}).

In section 5, we give the definition for an ideal to be closed under sums. We first observe that being closed under sums is necessary for an ideal to be a precovering. Then we see that this property is sufficient under certain conditions (Theorem~\ref{precovering}).

Finally we revise the definition of being generated by a set of homomorphisms. We see that if we allow infinite direct sums in the factorization of elements of an ideal, the results we have luckily still hold. We conclude with some immediate questions that still need to be answered.

\section{Preliminaries}

Let $R$ be an arbitrary ring and $R$-mod denote the category of left $R$ modules. Throughout the paper we will denote $\Ext^1_{R}$ as $\Ext$ and interpret $\Ext$ in terms of extensions of $R$ modules when necessary. An extension of $U$ by $V$ denoted by an element $\xi \in \Ext(U,V)$ is a short exact sequence,
$$\xymatrix{0 \ar[r] & U  \ar[r] & X \ar[r] & V \ar[r] & 0} $$
Two extensions are said to be equivalent if there is a homomorphism making the diagram,
$$\xymatrix{0 \ar[r] & U \ar@{=}[d] \ar[r] & X \ar[d] \ar[r] & V \ar@{=}[d] \ar[r] & 0\\
0 \ar[r] & U  \ar[r] & Y \ar[r] & V \ar[r] & 0} $$
commutative. Then one can easily see that there is an equivalence relation on $\Ext(U,V)$ for any given pair of modules and that $\Ext:R-mod \times R-mod \rightarrow Ab$ is a bifunctor covariant in the first component and contravariant in the second. Then $\Ext(f,g)$ will be calculated by using a pushout along $g$ followed by a pullback along $f$ (or equivalently a pullback along f followed by a pushout along g). Note that $\Ext(f,g)(\xi)=0$ means that the extension  $\Ext(f,g)(\xi)$ is split exact.

Throughout the paper an additive subfunctor $\ci$ of the $\Hom_{R}$ functor will be called an ideal. This is Herzog's definition first given in \cite{Herzog}. As a consequence one can easily observe that a class $\ci$ of $R$-homomorphisms will form an ideal if it satisfies the following conditions,
\begin{itemize}
\item If $f,g$ is in $\ci$ with the same domain and codomain then $f+g$ is in $\ci$.
\item If $g$ is in $\ci$ then for any $R$-homomorphisms $f,h$, $f \circ g \circ h$ is in $\ci$ (assuming the domains and codomains are suitable for the composition).
\end{itemize}

The definitions of a precover(right-approximation) and preenvelope(left-approximation) are carried over to the ideal case as given below.

 \begin{definition}Let $\ci \subseteq \Hom_{R}$ be an ideal and $M$ be a left R-module. An $\ci$-precover of $M$ is a morphism $i:I \rightarrow M$
 such that any $i': I' \rightarrow M$ from $I$ factors through $i$ as seen in the following diagram,
 $$\xymatrix{& I' \ar@{.>}[ld] \ar[d]^{i'} \\
 I \ar[r]^{i} & M} $$
 A $\cal{J}$-preenvelepe is defined similarly.
 \end{definition}

\begin{definition}Given two ideals $\ci, \cal{J} \subseteq \Hom_{R}$ of R-modules define,\\
 \begin{center}
 $ \cip$$= \{ j | \Ext^{1}(i,j)=0 \text{ for  all \,} i \in \ci \}$ \\
 $^{\perp} \cal{J}$$=\{i | \Ext^{1}(i,j)=0 \text{ for  all \,} j \in  \cal{J}  \}$
 \end{center}
 An ideal cotorsion pair in R-Mod is a pair $(\ci,\cal{J})$ of ideals such that $\ci^{\perp}=\cal{J}$ and $^{\perp}\cal{J}=\ci$.
 \end{definition}

\section{Properties of an ideal generated by a set}

In this section we observe how the elements of an ideal generated by a set can be factored through a certain kind of homomorphism. This observation helps us to identify the elements of $\cip$. We finish the section by Lemma~\ref{factor} which will be the main tool while proving $\cip$ to be preenveloping.

\begin{remark}Let $\ci$$=<f>$ where $f: M \rightarrow N$ then $\varphi:U \rightarrow V$ is in $I$ if and only if it has a factorization of the form,
$$\xymatrix{U \ar[r] &M^{m} \ar[r]^{f_{ji}} & N^{n} \ar[r] & V  }$$
for some $1 \leq m,n$ where $f_{ji}$ has entries either equal to $f$ or $0$.
\end{remark}

\proof Let $S=\{\varphi \, | \, \varphi \, \text{has the desired factorization property}\}$. Clearly $f \in \ci$ and $S \subseteq \ci$, so it is enough  to prove that $S$ is an ideal. Let $g,g' \in S$ where,
$$\xymatrix{g: U \ar[r]^{\alpha} &M^{m} \ar[r]^{g_{ji}} & N^{n} \ar[r]^{\beta} & V  } $$
and
$$\xymatrix{g': U \ar[r]^{\tilde{\alpha}} & M^{\tilde{m}} \ar[r]^{g'_{ji}}  & N^{\tilde{n}} \ar[r]^{\tilde{\beta}} & V  } $$
then $g+g'$ has the following factorization,
$$\xymatrix{U \ar[r]^-{(\alpha,\tilde{\alpha})} & M^{m+\tilde{m}} \ar[r]^{h}  & N^{n+\tilde{n}} \ar[r]^-{(\beta,\tilde{\beta})} & V  } $$
where
$$(\alpha,\tilde{\alpha})(u)=(\alpha(u),\tilde{\alpha}(u)),$$
$$(\beta,\tilde{\beta})(n,\tilde{n})=\beta(n)+\tilde{\beta}(n)$$
and
$$h=\begin{bmatrix}
(g_{ji}) & 0 \\
0  & (g'_{ji}) \\
\end{bmatrix}$$
for each $1 \leq j \leq n+\tilde{n}$ and $1 \leq i \leq m+\tilde{m}$. Given $u\in U$ say $\alpha(u)=(x_{1},...,x_{m})$ and $\tilde{\alpha}(u)=(y_{1},...,y_{\tilde{m}})$ then,
$$\xymatrix{u \ar@{|->}[r]^-{(\alpha,\tilde{\alpha})} & (\alpha(u),\tilde{\alpha(u)}) \ar[r]^-{h} & (g_{ji}(x_{j}),g'_{ji}(y_{j})) \ar[r]^-{(\beta,\tilde{\beta})} &\beta(g_{ji}(x_{j}))+ \tilde{\beta}(g'_{ji}(y_{j})) }$$
where notice that,
$$g(u)+ g'(u)= \beta(g_{ji}(x_{j}))+ \tilde{\beta}(g'_{ji}(y_{j})) $$
Hence we conclude that we have a factorization of $g+g'$. \hfill $\Box$

\begin{remark}Let $\ci$$=<f^{k}>_{k \in K}$ be generated by a set of homomorphisms where $f^{k}:M_{k}\rightarrow N_{k}$ then $\varphi:U \rightarrow V$ is in $\ci$ if and only if it has a factorization as following,
$$\xymatrix{U \ar[r] & M_{k_{1}}^{m_{1}}\oplus ...\oplus M_{k_{t}}^{m_{t}}  \ar[r]^{(h_{ji})} & N_{k_{1}}^{n_{1}}\times ...\times N_{k_{t}}^{n_{t}} \ar[r] &V }
$$ where $k_{1},...,k_{t} \in K$ and

\begin{align*}
h_{ji}&=f^{k_{1}} \text{\, or \,} 0 \text{\, for \,}  1 \leq j \leq n_{1}, 1\leq i \leq m_{1} ,\\
\hspace*{0.5in} . \\
\hspace*{0.5in} .\\
\hspace*{0.5in} . \\
h_{ji}&=f^{k_{t}} \text{\,or \,} 0 \text{\, for \,} (n_{1}+...+n_{t-1})\leq j \leq(n_{1}+...+n_{t}), \, (m_{1}+...+m_{t-1})\leq i \leq ( m_{1}+...+m_{t})
\end{align*}
and $h_{ji}$ can be viewed as a matrix with entries,
$$ \begin{bmatrix}
       f_{ji}^{k_{1}} & 0 & \cdots  & 0    \\
       0 & f_{ji}^{k_{2}} & \cdots  & 0 \\
       \vdots & \vdots  & \ddots & \vdots \\
       0    & 0        & \cdots & f_{ji}^{k_{t}}
     \end{bmatrix}
$$
 \label{setgenerated}
\end{remark}

\proof Very similar to the case where $\ci$ is generated by a single homomorphism.

\begin{lemma}Let $\ci$$=< f_{k} \, | \, k \in K >$ and $f: M_{k} \rightarrow N_{k}$, $g: U \rightarrow V$ be homomorphisms of R-modules then the following are equivalent,
\begin{enumerate}
 \item $g$ is in $\cip$.
 \item Given any s.e.s.,
$$\xymatrix{ 0 \ar[r] & U \ar[r] & X \ar[r] & Y \ar[r] & 0} $$
the s.e.s. we get by using the pushout along $g$,
$$\xymatrix{ 0 \ar[r] & U \ar[r] \ar[d]^{g} & X \ar[r]  \ar[d] & Y \ar[r] \ar@{=}[d] & 0\\
 0 \ar[r] & V \ar[r] & Q \ar[r] & Y \ar[r] & 0} $$
can be completed to a commutative diagram for any $N_{k} \rightarrow Y$ where $k \in K$ as shown below,
$$\xymatrix{ &  &  & M_{k} \ar@{.>}^{\displaystyle{\circlearrowleft}}[ddl] \ar[d]^{f_{k}} & \\
  &  &  & N_{k} \ar[d] & \\
 0 \ar[r] & V \ar[r] & Q \ar[r] & Y \ar[r] & 0} $$
 \item Let $0\rightarrow V\rightarrow Q\rightarrow Q/V\rightarrow 0$ be the s.e.s. obtained by using the pushout along $g$, then for any $k \in K$ the composition $ \Hom(N_{k}, Q/V) \rightarrow \Hom(M_{k} , Q/V) \rightarrow  \Ext^{1}(M_{k},V)$
 (obtained from $\xymatrix{ M_{k} \ar[r]^{f_{k}}& N_{k} \ar[r]& Q/V}$) is the zero map for any homomorphism $N_{k}\rightarrow Q/V$.
 \end{enumerate}
\label{iperp}
\end{lemma}

\proof ($1 \Rightarrow 2$) Assume that the first statement is correct. Since $g \in \cip$, then $\Ext^{1}(\tilde{f},g)(\xi)=0$ for any $\tilde{f} \in \ci$. That is the resulting s.e.s we get by computing the pushout along $g$ followed by the pullback along $\tilde{f}$ is split exact. But since $M_{k} \rightarrow N_{k} \rightarrow Y$ is in $\ci$,
$$\xymatrix{ 0 \ar[r] & V \ar[r] \ar@{=}[dd] & P \ar[dd] \ar[r]& M_{k} \ar@/_/[l] \ar@{.>}[ddl] \ar[d]^{f_{k}} \ar[r]& 0\\
  &  &  & N_{k} \ar[d] & \\
 0 \ar[r] & V \ar[r] & Q \ar[r] & Y \ar[r] & 0} $$
we get a split exact sequence on the upper row. Hence we can complete the triangle with the dotted map (which is obtained by the composition $M_{k} \rightarrow P \rightarrow Q$) above.

($2 \Rightarrow 1$) Assume now that the second property holds for $g: U \rightarrow V$. We need to prove that $\Ext^{1}(\tilde{f},g)(\xi)=0$ where $$\xymatrix{ \tilde{f}: M_{j} \oplus M_{i} \ar[r]^-f & N_{j} \oplus N_{i} \ar[r] & Y }$$
in $\ci$ and for any s.e.s. $\xi: 0 \rightarrow U \rightarrow X \rightarrow Y \rightarrow 0$ (then the proof follows very similarly for an arbitrary $\tilde{f} \in \ci$). By the previous proposition, $f=\left(
                                                        \begin{array}{cc}
                                                          f_{j} & 0 \\
                                                          0& f_{i}\\
                                                        \end{array}
                                                      \right)$
where $f_{j},f_{i}$ are from the set of generators of $\ci$. Given any $\varphi: N_{j}\oplus N_{i} \rightarrow Y$, we define $\varphi_{j}: N_{j}\rightarrow Y $ and $\varphi_{i}: N_{i}\rightarrow Y $ such that $\varphi_{j}$ is the restriction of $\varphi$ to $N_{j} \oplus 0$ and similarly $\varphi_{i}$ is the restriction of $\varphi$ to $0 \oplus N_{i}$. Then by assumption there exists $\alpha_{j}$ and $\alpha_{i}$ making the following diagrams commutative,
$$\xymatrix{ &  &  & M_{j} \ar@{.>}^{\displaystyle{\circlearrowleft}}[ddl]_{\alpha_{j}} \ar[d]^{f_{j}} & \\
  &  &  & N_{j} \ar[d]^{\varphi_{j}} & \\
 0 \ar[r] & V \ar[r] & Q \ar[r] & Y \ar[r] & 0} $$
 and
$$\xymatrix{ &  &  & M_{i} \ar@{.>}^{\displaystyle{\circlearrowleft}}[ddl]_{\alpha_{i}} \ar[d]^{f_{i}} & \\
  &  &  & N_{i} \ar[d]^{\varphi_{i}} & \\
 0 \ar[r] & V \ar[r] & Q \ar[r] & Y \ar[r] & 0} $$
Then it is easy to see that the map $\alpha: M_{j} \oplus M_{i} \rightarrow Q $ defined as $\alpha(x_{j},x_{i})=\alpha_{j}(x_{j})+\alpha_{i}(x_{i})$ makes the following diagram commutative,
$$\xymatrix{ &  &  & M_{j}\oplus M_{i} \ar[ddl]_{\alpha} \ar[d]^{f} & \\
  &  &  & N_{j}\oplus N_{i} \ar[d]^{\varphi} & \\
 0 \ar[r] & V \ar[r] & Q \ar[r] & Y \ar[r] & 0} $$
If we compute the pullback along $\tilde{f}$ we get the following diagram,
$$\xymatrix{ P \ar[r] \ar[d] &  M_{j}\oplus M_{i}  \ar[d]^{\tilde{f}} \\
Q \ar[r] & Y } $$
Now using the commutativity of the previous diagram, we get the following commutative diagram,
$$\xymatrix{  M_{j}\oplus M_{i} \ar@/^/^{id}_{\circlearrowleft}[drr]  \ar@/_/_{\alpha}^{\circlearrowright}[ddr] \ar@{.>}[dr]^{\psi} & &\\
&P \ar[r] \ar[d] &  M_{j}\oplus M_{i}  \ar[d] \\
&Q \ar[r] & Y } $$
So by the universal property of pushout diagrams we conclude that there exists a homomorphism $\psi$ such that $\xymatrix{ M_{j}\oplus M_{i}  \ar^-\psi[r] & P \ar[r] &  M_{j}\oplus M_{i} } $ is the identity morphism. Hence looking at the s.e.s. obtained by the pullback along $\tilde{f}$,
$$\xymatrix{0 \ar[r] & U \ar[r] \ar@{=}[d] & P \ar[r] \ar[d] &  M_{j}\oplus M_{i} \ar[r] \ar^{\tilde{f}}[d] & 0 \\
0 \ar[r] & U \ar[r] & Q \ar[r] & Y \ar[r]  & 0} $$
we conclude that the upper row is split exact. That is $\Ext^{1}(\tilde{f},g)$ maps $\xi$ to a split exact sequence, i.e. $\Ext^{1}(\tilde{f},g)(\xi)=0$ for any s.e.s. $\xi$.

($2 \Leftrightarrow 3$) Assume that the second property holds. Given any s.e.s. $\xi: 0 \rightarrow U \rightarrow X \rightarrow Y \rightarrow 0$  by using the pushout along $g$ we get,
$$\xymatrix{ 0 \ar[r] & U \ar[r] \ar[d]^{g} & X \ar[r]  \ar[d] & Y \ar[r] \ar@{=}[d] & 0\\
 0 \ar[r] & V \ar[r] & Q \ar[r] & Y\cong Q/V \ar[r] & 0} $$
by assumption the lower row can be completed to a commutative diagram for any $N_{k} \rightarrow  Q/V$ as shown below,
$$\xymatrix{ &  &  & M_{k} \ar@{.>}^{\displaystyle{\circlearrowleft}}[ddl] \ar[d]^{f_{k}} & \\
  &  &  & N_{k} \ar[d] & \\
 0 \ar[r] & V \ar[r] & Q \ar[r] & Q/V \ar[r] & 0} $$
So we get the diagram,
$$\xymatrix{ &\Hom(N_{k}, Q/V) \ar[d] & \\
\Hom(M_{k},Q)\ar[r] &\Hom(M_{k},Q/V) \ar[r] & \Ext^{1}(M_{k},V)} $$
with an exact row. But now our assumption holds if and only if the following composition,
$$\xymatrix{\Hom(N_{k},Q/V)\ar[r] &\Hom(M_{k},Q/V) \ar[r] & \Ext^{1}(M_{k},V) } $$
is the zero map. \hfill $\Box$

\begin{corollary}Let $\ci$ be as in lemma~\ref{iperp} and $g: U \rightarrow V$ be in $\cip$. If $V^{'} \subseteq V$ is a submodule such that $g(U)\subseteq V^{'} \subseteq V$ and the map $\Ext^{1}(M_{k},V^{'}) \rightarrow \Ext^{1}(M_{k},V)$ is an injection for any $k \in K$ then $g: U \rightarrow V^{'} $ is in $\cip$ as well.\label{extinjection}
\end{corollary}

\proof Notice that by lemma~\ref{iperp} we conclude that $g:U \rightarrow V' \subseteq V$ is in $\cip$ if and only if the following composition is $0$ for any given $N_{k} \rightarrow Q^{'}/V^{'}$ and any given $k \in K$,
$$\xymatrix{\Hom(N_{k},Q^{'}/V^{'})\ar[r] &\Hom(M_{k},Q^{'}/V^{'}) \ar[r] & \Ext^{1}(M_{k},V^{'}) } $$
which is induced from following s.e.s.,
$$0\rightarrow V^{'}\rightarrow Q^{'}\rightarrow Q^{'}/V^{'}\rightarrow 0$$
Then the following diagram,
$$\xymatrix{ 0 \ar[r] & U \ar[r] \ar[d]^{g} & E(U)\ar[r]  \ar[d] & E(U)/U \ar[r] \ar[d]^{\cong} & 0\\
 0 \ar[r] & V^{'} \ar[r] \ar@{^{(}->}[d] & Q^{'} \ar[r] \ar[d] &  Q^{'}/V^{'} \ar[r] \ar[d]^{\cong} & 0\\
 0 \ar[r] & V \ar[r] & Q \ar[r] & Q/V \ar[r] & 0} $$
gives us,
$$\xymatrix{\Hom(N_{k},Q'/V^{'})\ar[r] \ar[d]^{\cong} &\Hom(M_{k},Q'/V^{'}) \ar[r] \ar[d]^{\cong} & \Ext^{1}(M_{k},V^{'}) \ar[d] \\
\Hom(N_{k},Q/V)\ar[r] &\Hom(M_{k},Q/V) \ar[r] & \Ext^{1}(M_{k},V)} $$
We notice that if $\Ext^{1}(M_{k},V^{'}) \rightarrow \Ext^{1}(M_{k},V)$ is an injection for every $k \in K$ then the composition on the top row is $0$ for every $k \in K$. By lemma~\ref{iperp} we conclude that $g:U \rightarrow V'$ is in $\cip$. \hfill $\Box$

\begin{lemma}Let $\ci$ be as in lemma~\ref{iperp}. If $g_{i}: U \rightarrow V_{i}$, $i \in I$ are each in $\cip$ then $g: U \rightarrow \prod\limits_{i \in I} V_{i}$ where $g(x)=(g_{i}(x))_{i \in I}$ is in $\cip$.\label{prod}
\end{lemma}

\proof Assume $g_{i}: U \rightarrow V_{i}$, $i \in I$ are each in $\cip$. That is $\varphi_{i}=\Ext^{1}(\tilde{f},g_{i})$ is the zero map for any $\tilde{f}:Z \rightarrow Y $ in $\ci$. Since $\Ext^{1}(Z,\prod\limits_{i \in I} V_{i})\cong \prod\limits_{i \in I}\Ext^{1}(Z,V_{i})$ we have the following diagram,
$$\xymatrix{ & \Ext^{1}(Y,U) \ar^{\varphi_{j}}[ddr] \ar_{\varphi_{i}}[ddl] \ar@{.>}[d]& \\
& \prod\limits_{i \in I}\Ext^{1}(Z,V_{i}) \ar_{\pi_{j}}[dr] \ar^{\pi_{i}}[dl] & \\
\Ext^{1}(Z,V_{i})  & & \Ext^{1}(Z,V_{j})} $$
So there exists $\varphi: \Ext^{1}(Y,U) \rightarrow \prod\limits_{i \in I}\Ext^{1}(Z,V_{i}) $ which makes the above diagram commutative. Now given any $\eta \in \Ext^{1}(Y,U)$ and say $\varphi(\eta)=(\xi_{i})_{i \in I}$ then,
$$0=\varphi_{j}(\eta)=\pi_{j}(\varphi(\eta))=\pi_{j}((\xi_{i})_{i \in I})=\xi_{j} $$
That is $\eta=0$, which gives us that $\varphi$ or $\Ext(\tilde{f},g)$ is the zero map. Hence $g \in \cip$. \hfill $\Box$

\begin{lemma}Let $\ci$ be as in lemma~\ref{iperp} and $g: U \rightarrow V$ in $\cip$. Then $g$ can be factored through $V^{'}$ such that,
$$\xymatrix{U \ar[r] \ar@/^1pc/^{g}[rr] &V^{'} \ar@{^{(}->}[r] & V} $$
where the cardinality of $V^{'}$ is bounded by a cardinal number $\kappa$ which depends only on $|U|$ and $\ci$. \label{factor}
\end{lemma}

\proof First we need to show that $g$ is in $\cip$ if and only if $\Ext(\tilde{f},g)(\xi^{'})=0$ for the short exact sequence $\xi^{'}: 0 \rightarrow U \rightarrow E(U) \rightarrow E(U)/U \rightarrow 0 $ where $E(U)$ is the injective hull of $U$. One way is obvious. To show the other way let,
$$\xi: 0 \rightarrow U \rightarrow X \rightarrow Y \rightarrow 0$$
be any short exact sequence. Since $E(U)$ is injective we get the following commutative diagram,
$$\xymatrix{ 0 \ar[r] & U \ar[r] \ar@{=}[d] & X \ar[r]  \ar[d] & Y \ar[r] \ar@{.>}^{\exists k}[d] & 0\\
 0 \ar[r] & U \ar[r] & E(U) \ar[r] & E(U)/U \ar[r] & 0} $$
where $k: Y \rightarrow E(U)/U$ is induced from $X \rightarrow E(U)$.
Now we look at $\Ext^{1}(k,id_{U})(\xi^{'})$,
$$\xymatrix{ 0 \ar[rr] && U \ar[rr] \ar@{=}[dd] && X \ar[rr] \ar[dd] && Y \ar[r] \ar[dd]^{k} & 0 && \\
& 0 \ar[rr] && U \ar[rr] \ar[dl]^{id_{U}} && E(U) \ar[rr] \ar[dl] && E(U)/U \ar[r] \ar[dl] & 0 \\
 0 \ar[rr] && U \ar[rr] && E(U) \ar[rr] && E(U)/U \ar[r]  & 0 &&}$$
That is $\Ext^{1}(k,id_{U})(\xi^{'})=\xi$. Note that the upper row is a pullback along $k$ since $X \rightarrow Y$ and $E(U) \rightarrow E(U)/U$ are epimorphisms and $\xymatrix{U \ar[r]^{id_{U}} & U}$ is an isomorphism.
Then,
$$\Ext^{1}(\tilde{f},g)(\xi)= \Ext^{1}(\tilde{f},g) \circ \Ext^{1}(k,id_{U})(\xi^{'}) = \Ext^{1}(k \circ \tilde{f},g)(\xi^{'})=0 $$
Hence $g \in \cip$.\\
Now we want to construct a ``small'' enough $V^{'}$ such that $g: U \rightarrow V^{'} \subset V$ is in $\cip$. By corollary~\ref{extinjection} it is enough to show that $\Ext^{1}(M_{k},V^{'}) \rightarrow \Ext^{1}(M_{k},V)$ is an injection for any $k\in K$. But notice that this holds if the following map is an injection,
$$\prod\limits_{k \in K} \Ext^{1}(M_{k},V^{'}) \rightarrow \prod\limits_{k \in K} \Ext^{1}(M_{k},V) $$
But we have the following commutative diagram,
$$\xymatrix{ \Ext^{1}(\bigoplus\limits_{k\in K}M_{k},V^{'}) \ar[r]^{\cong} \ar[d] & \prod\limits_{k \in K} \Ext^{1}(M_{k},V^{'}) \ar[d] \\
\Ext^{1}(\bigoplus\limits_{k\in K}M_{k},V) \ar[r]^{\cong} & \prod\limits_{k \in K} \Ext^{1}(M_{k},V)}$$
So we conclude that to show $g:U \rightarrow V^{'}$ in $\cip$ it is enough to show the left column is injective.

Now we will construct the desired $V^{'}$. Given $i\in K$. Let $0 \rightarrow K_{i} \rightarrow P_{i} \rightarrow M_{i} \rightarrow 0$ be a partial projective resolution of $M_{i}$. Then we obtain the following partial projective resolution for $\bigoplus\limits_{i \in K} M_{i} $,
$$\xymatrix{0 \ar[r] & K=\bigoplus\limits_{i \in K} K_{i} \ar[r] & P=\bigoplus\limits_{i \in K} P_{i} \ar[r] & \bigoplus\limits_{i \in K} M_{i} \ar[r] &0 } $$
We will construct an ascending chain of modules to obtain such a ``small'' $V^{'}$. Let $g(U)=V_{0}$ and for every $K \rightarrow V_{0}$ morphism that has an extension $P \rightarrow V$, choose one such extension $\alpha: P \rightarrow V$. So we have the following diagram,
$$\xymatrix{P \ar[rd] & \\
K \ar@{^{(}->}[u] \ar[r]&  V_{0}\subset V } $$
Define $V_{1}= \sum \alpha(P)$ where the sum is over all such chosen extensions $\alpha: P \rightarrow V$ for each $K \rightarrow V_{0} $. Then $V_{0} \subset V_{1}$. Now we construct $V_{2}$ in a similar way. So we get an ascending chain of modules,
$$V_{0} \subset V_{1}  \subset ... \subset V_{\omega} \subset V_{\omega+1} \subset ...\subset V_{\beta} $$
where $\beta$ is the least cardinal number with $|K|< \beta$ . Define $V_{\lambda}=\bigcup\limits_{\alpha < \lambda} V_{\alpha}$ if $\lambda \leq \beta$ is a limit ordinal and let $V^{'}=V_{\beta}$.

Now we will show that $g: U \rightarrow V^{'}$ is in $\cip$. By the previous observation all we need is,
$$\Ext^{1}(\bigoplus\limits_{k\in K}M_{k},V^{'}) \rightarrow \Ext^{1}(\bigoplus\limits_{k\in K}M_{k},V)$$
to be injective. Given a morphism $\xymatrix{K \ar[r]^-{\varphi} & V^{'} \subset V}$ that has an extension $\xymatrix{P \ar[r]^{\Phi} & V} $ we want to show that then there is an extension $P\rightarrow V^{'}$. But now since $|K| < \beta$ we conclude that $Im(\varphi)\subseteq V_{\alpha}$ for some $\alpha$ such that $|\alpha| < \beta$ hence by the construction of the ascending chain we can extend $\varphi$ to $P \rightarrow V_{\alpha+1} \subset V^{'}$. This shows that
$$\Ext^{1}(\bigoplus\limits_{k\in K}M_{k},V^{'}) \rightarrow \Ext^{1}(\bigoplus\limits_{k\in K}M_{k},V)$$
is an injection. Now by corollary~\ref{extinjection} we conclude that $g: U \rightarrow V^{'} \subset V$ is in $\cip$. Moreover note that the cardinality of $V_{\alpha+1}= \sum \Phi(P)$ (where sum is over all chosen extensions $\Phi: P \rightarrow V$ for each $K \rightarrow V_{\alpha}$) is bounded by,
$$|V_{\alpha+1}| \leq |P|^{|V_{\alpha}|^{|K|}  }$$
Since,
$$|\Hom(K,V_{\alpha}) | \leq |V_{\alpha}^{K}|=|V_{\alpha}|^{K} $$
and
$$|\Phi(P)| \leq |P| $$
We find a bound on $|V^{'}| \leq \sum\limits_{\alpha < \beta} |V_{\alpha}|$. So we conclude that for any given $U$, if $g: U \rightarrow V$ is in  $\cip$ we can find a factorization $ U \rightarrow V^{'} \rightarrow V$ where $|V^{'}| \leq \kappa$ for some cardinal number $\kappa$ that depends on $|U|$, $|K|$ and $|P|$ and so only on $|U|$ and $\ci$. \hfill $\Box$

\section{Main result}

We now prove the main theorem. This result was motivated by the theorem of Eklof-Trlifaj's (Thm.10, \cite{Eklof}).
\begin{theorem} If $\ci$ is generated by a set then $\cip$ is a preenveloping class.\label{preenveloping}\end{theorem}

\proof Given a R-module $U$ by lemma~\ref{factor} we find a cardinal number $\kappa$ with the desired properties given in the lemma . We will use a similar argument to Rada-Saor\'{i}n's from \cite{RadaSaorin}. Let $\{ g_{j}\}_{j \in J}$ be the set of all the homomorphisms $g_{j}: U \rightarrow V_{j}$ in $\cip$ (up to isomorphism) where $|V_{j}| \leq \kappa$. Then by lemma~\ref{prod} the homomorphism,
$$\xymatrix{U \ar[r]^{\prod g_{j}} & \prod\limits_{j \in J} V_{j} }$$
is in $\cip$, moreover we claim that it is an $\cip$-preenvelope of $U$.\\
Given any $g: U \rightarrow V$ in $\cip$ by lemma~\ref{factor} we get a factorization,
$$\xymatrix{U \ar[r] \ar@/^1pc/[rr]^{g} & V' \ar@{^{(}->}[r]& V} $$
where $|V'| \leq \kappa$. Then $g: U  \rightarrow V' \subset V$ is isomorphic to $g_{j}$ for some $j$. That is there exists a map making the following commutative,
$$\xymatrix{ U  \ar[r]^{g_{j}} \ar[d]_{g} & V_{j} \ar[dl]_{\displaystyle\circlearrowright}  \\
V' & } $$
Now we get the following commutative diagram,
$$\xymatrix{ U  \ar[r]^{g_{j}} \ar[d] & \prod\limits_{j \in J} V_{j} \ar[dl] \\
V_{j} \ar[d]& \\
V' \ar[d]& \\
V & } $$
We conclude that the following diagram is commutative where we use the composition of maps $\prod\limits_{j \in J} V_{j} \rightarrow V_{j} \rightarrow V' \rightarrow V$ ,
$$\xymatrix{ U  \ar[r]^{g_{j}} \ar[d]_{g} & \prod\limits_{j \in J} V_{j} \ar[dl]_{\displaystyle\circlearrowright} \\
V & } $$
That is $\cip$ is a preenveloping class. \hfill $\Box$

\section{A necessary condition for $\ci$ to be precovering}

\begin{definition} The ideal $\ci$ is said to be closed under sums if it satisfies one of the following equivalent conditions,
 \begin{itemize}
 \item If $(f_{j})_{j \in J}$, $f_{j}: M_{j} \rightarrow N_{j}$ is any family of elements of $\ci$ then $\bigoplus\limits_{j \in J}f_{j}: \bigoplus\limits_{\j \in J} M_{j} \rightarrow \bigoplus\limits_{\j \in J} N_{j}$ is in $\ci$.
\item If $(g_{j})_{j \in J}$, $g_{j}: M_{j} \rightarrow N$ is any family of elements of $\ci$ then $ g: \bigoplus\limits_{\j \in J} M_{j} \rightarrow N$ defined by $g((x_{j})_{j \in J})=\sum\limits_{j \in J}g_{j}(x_{j})$ is in $\ci$.
\end{itemize}
\end{definition}

\begin{theorem} If $\ci$ is the closure under direct sums of the ideal generated by a single homomorphism $f:M \rightarrow N$ then $\ci$ is a precovering ideal.\label{single}\end{theorem}

\proof Given an arbitrary $R$-module $V$, we consider the homomorphism $\alpha: M^{(\Hom(N,V))} \rightarrow V $ defined by,
$$\alpha((x_{g})_{g\in \Hom(N,V)})=\displaystyle\sum\limits_{g\in \Hom(N,V)}  g(f(x_{g})) $$
which is in $\ci$ since it is closed under sums. Moreover we claim that it is a $\ci$-precover of $V$. Given a homomorphism,
$$\xymatrix{\bigoplus\limits_{j \in J} M \ar[r] & \bigoplus\limits_{j \in J}N \ar[r]^{h} & V }$$
 in $\ci$. Define $h_{j}: N \rightarrow V $ such that $h=\sum\limits_{j \in J} h_{j}$  and $\beta_{j}:M \rightarrow M^{(\Hom(N,V))} $ such that $\beta_{j}(x)$ is the element of $M^{(\Hom(N,V))}$ whose all entries are $0$, except the one that corresponds to $h_{j}$, which is $x$. With the maps defined above following diagram commutes,
$$\xymatrix{ & M \ar[ddl]_{\beta_{j}} \ar[d]^{f} \\
& N \ar[d]^{h_{j}} \\
M^{(\Hom(N,V))} \ar[r]^-\alpha &  V } $$

Now we define $\beta$,
$$\xymatrix{ & \bigoplus\limits_{j \in J} M \ar@{.>}[ddl]_{\beta} \ar[d]^-{\bigoplus f} \\
& \bigoplus\limits_{j \in J}N \ar[d]^{h} \\
M^{(\Hom(N,V))} \ar[r]^-\alpha &  V } $$
where $\beta((x_{j})_{j \in J})=\sum\limits_{j\in J} \beta_{j}((x_{j})$. Then notice that,
\begin{align*}
\alpha(\beta((x_{j})_{j \in J})) &= \alpha(\sum\limits_{j\in J} \beta_{j}(x_{j}))=\sum\limits_{j\in J} \alpha(\beta_{j}(x_{j}))  \\
&=\sum\limits_{j\in J} h_{j}(f(x_{j}))= h((f(x_{j}))_{j\in J}) \\
&=(h ( \oplus f(x_{j})_{j \in J}))
\end{align*}
So the above diagram is commutative and we conclude that $V$ has an $\ci$-precover. That is $\ci$ is a precovering ideal. \hfill $\Box$

\begin{theorem} If $\ci$ is the closure under direct sums of the ideal generated by a set then $\ci$ is a precovering ideal.\label{precovering} \end{theorem}

\proof The proof follows very similarly to that of theorem~\ref{single}. \hfill $\Box$

\begin{proposition}Let $\ci=$$< f_{s} >_{s\in S}$ generated by a set and $\ci^{'}$ be the smallest ideal that contains $\ci$ and closed under sums. Then $\cip= (\ci^{'})^{ \perp}$.\label{extendedperp}
\end{proposition}

\proof One way of the inclusion is easy. By definition $\ci \subseteq \ci^{'}$ implies $(\ci^{'})^{ \perp} \subseteq \cip$.

To prove the other way of inclusion let $g \in \cip$ where $g: U \rightarrow V$. Then notice that $\Ext^{1}(f_{s},g)=0$ for all $f_{s}$. Now given any $T \subseteq S$ and a homomorphism,
$$\xymatrix{\bigoplus M_{t} \ar[r]^{\oplus f_{t}} &  \bigoplus N_{t} }$$
We want to prove that,
$$\xymatrix{\Ext^{1}(\oplus f_{t}, g): \Ext^{1}(\oplus N_{t},U) \ar[r] & \Ext^{1}(\oplus M_{t},V)}$$
is the zero map. We observe the commutative diagram,
$$\xymatrix{ \Ext^{1}(\bigoplus\limits_{t\in T}N_{t},U) \ar[r]^{\cong} \ar[d] & \prod\limits_{t\in T} \Ext^{1}(N_{t},U) \ar[d] \\
\Ext^{1}(\bigoplus\limits_{t\in T}M_{t},V)  \ar[r]^{\cong} & \prod\limits_{t\in T} \Ext^{1}(M_{t},V) }$$
Notice that the right column is the zero map since $\Ext^{1}(f_{t},g)=0$ for all $t \in T$. Hence we conclude that $\Ext^{1}(\oplus f_{t}, g)=0$. Now given an arbitrary homomorphism in $\ci^{'}$,
$$\xymatrix{U \ar[r]^-{k} & \bigoplus M_{t} \ar[r]^-{\oplus f_{t}} & \bigoplus N_{t} \ar[r]^{h} & V } $$
where $t \in T \subseteq S$. But we have,
$$\Ext^{1}(h\circ \oplus f_{t} \circ k,g)=\Ext^{1}(k,id)\circ \Ext^{1}(\oplus f_{t} ,g) \circ \Ext^{1}(h,id)=0$$
Hence we conclude that $\cip \subseteq (\ci^{'})^{ \perp}$. \hfill $\Box$

\section{Ideals generated by a set in the extended sense}

We revise our definition of $\ci$ being generated by a set.

\begin{definition} Let $(f_{s})_{s \in S}$ be a set of homomorphisms where $f_{s}: M_{s} \rightarrow N_{s}$. $\ci$ is said to be generated by $(f_{s})_{s \in S}$ in the extended sense if every $\tilde{f}:U \rightarrow V$  in $\ci$ has a factorization,
$$\xymatrix{U \ar[r] & \bigoplus\limits_{s\in S}M_{s}^{\kappa_{s}} \ar[r]^{\bigoplus f_{s}}  & \bigoplus\limits_{s\in S}N_{s}^{\kappa_{s}} \ar[r] & V } $$
where $\kappa_{s}$ is a cardinal number for each $s \in S$.
\end{definition}

\begin{remark} If $\ci$ is generated by a set in the extended sense, then it is closed under sums.\label{closedsum}
\end{remark}

\begin{corollary}Let $\ci$ be generated by a set of homomorphisms in the extended sense, then $\cip$ is a preenveloping ideal.
\end{corollary}

\proof Theorem~\ref{preenveloping} and proposition~\ref{extendedperp} gives us the result. \hfill $\Box$

\begin{corollary}Let $\ci$ be generated by a set of homomorphisms in the extended sense, then $\ci$ is a precovering ideal.
\end{corollary}

\proof Theorem~\ref{precovering} and remark~\ref{closedsum} gives us the result. \hfill $\Box$

There are questions still need to be answered when $\ci$ is an ideal generated by a set in the extended sense, such as whether $(\ci,\cip)$ is a cotorsion ideal pair (i.e is $\ci=$$^{\perp}$$(\cip)$) and what the necessary conditions are for completeness if $(\ci,\cip)$ is a cotorsion ideal pair.

\textbf{Acknowledgements.} The author is very thankful to her PhD. advisor Prof. Edgar Enochs for his contributions and suggestions during the course of preparing this manuscript.

\end{document}